\pdfoutput=1
\RequirePackage{ifpdf}
\ifpdf 
\documentclass[pdftex]{sigma}
\else
\documentclass{sigma}
\fi

\newcommand*\pfqskip{8mu}
\catcode`,\active
\newcommand*\pfq{\begingroup
        \catcode`\,\active
        \def ,{\mskip\pfqskip\relax}%
        \dopfq
}
\catcode`\,12
\def\dopfq#1#2#3#4#5{%
        {}_{#1}\phi_{#2}\bigg(\genfrac..{0pt}{}{#3}{#4}\Big\rvert#5\biggr)%
        \endgroup
}
\newcommand*\pFqskip{8mu}
\catcode`,\active
\newcommand*\pFq{\begingroup
        \catcode`\,\active
        \def ,{\mskip\pFqskip\relax}%
        \dopFq
}
\catcode`\,12
\def\dopFq#1#2#3#4#5{%
        {}_{#1}F_{#2}\biggl[\genfrac..{0pt}{}{#3}{#4};#5\biggr]%
        \endgroup
}

\newtheorem*{Corollary}{Corollary}

\numberwithin{equation}{section}

\begin{document}

\allowdisplaybreaks

\renewcommand{\PaperNumber}{018}

\FirstPageHeading

\ShortArticleName{Bispectrality of the Complementary Bannai--Ito Polynomials}

\ArticleName{Bispectrality of the Complementary
\\
Bannai--Ito Polynomials}

\Author{Vincent X.~GENEST~$^\dag$, Luc VINET~$^\dag$ and Alexei ZHEDANOV~$^\ddag$}

\AuthorNameForHeading{V.X.~Genest, L.~Vinet and A.~Zhedanov}

\Address{$^\dag$~Centre de Recherches Math\'ematiques, Universit\'e de Montr\'eal,
\\
\hphantom{$^\dag$}~C.P.~6128, Succursale Centre-ville, Montr\'eal, Qu\'ebec, Canada, H3C 3J7}

\EmailD{\href{mailto:genestvi@crm.umontreal.ca}{genestvi@crm.umontreal.ca},
\href{mailto:luc.vinet@umontreal.ca}{luc.vinet@umontreal.ca}}

\Address{$^\ddag$~Donetsk Institute for Physics and Technology, Ukraine}
\EmailD{\href{mailto:zhedanov@yahoo.com}{zhedanov@yahoo.com}}

\ArticleDates{Received November 13, 2012, in f\/inal form February 27, 2013; Published online March 02, 2013}

\Abstract{A one-parameter family of operators that have the complementary Bannai--Ito (CBI)
polynomials as eigenfunctions is obtained.
The CBI polynomials are the kernel partners of the Bannai--Ito polynomials and also correspond to
a~$q\rightarrow-1$ limit of the Askey--Wilson polynomials.
The eigenvalue equations for the CBI polynomials are found to involve second order Dunkl shift
operators with ref\/lections and exhibit quadratic spectra.
The algebra associated to the CBI polynomials is given and seen to be a~deformation of the
Askey--Wilson algebra with an involution.
The relation between the CBI polynomials and the recently discovered dual $-1$ Hahn
and para-Krawtchouk polynomials, as well as their relation with the symmetric Hahn polynomials, is
also discussed.}

\Keywords{Bannai--Ito polynomials; quadratic algebras; Dunkl operators}

\Classification{33C02; 16G02}

\section{Introduction}
One of the recent advances in the theory of orthogonal polynomials (OPs) has been the discovery of
several new families of ``classical'' OPs that correspond to $q\rightarrow-1$ limits of
$q$-polynomials of the Askey scheme~\cite{VZhedanov-2011,Vinet-2012,Vinet-2011-2,Zhedanov-2011}.
The word ``classical'' here refers to the fact that in addition to obeying the three-term relation
\begin{gather*}
\mathcal{P}_{n+1}(x)+\beta_n\mathcal{P}_{n}(x)+\gamma_n\mathcal{P}_{n-1}(x)=x\mathcal{P}_{n}(x),
\end{gather*}
the polynomials $\mathcal{P}_n(x)$ also satisfy an eigenvalue equation of the form
\begin{gather*}
\mathcal{L}\mathcal{P}_{n}(x)=\lambda_n\mathcal{P}_n(x).
\end{gather*}
The novelty of these families of $-1$ orthogonal polynomials lies in the fact that for each family
the operator $\mathcal{L}$ is a~dif\/ferential or dif\/ference operator that also contains the
ref\/lection operator $Rf(x)=f(-x)$~\cite{Vinet-2011-3}.
Such dif\/ferential/dif\/ference operators are said to be of Dunkl type~\cite{Dunkl-1991},
notwithstanding the fact that the operators $\mathcal{L}$ dif\/fer from the standard Dunkl
operators in that they preserve the linear space of polynomials of any given maximal degree.
In this connection, these $-1$ OPs have also been referred to as Dunkl orthogonal polynomials.

With the discovery and characterization of these Dunkl polynomials, a~$-1$ scheme of OPs,
completing the Askey scheme, is beginning to emerge.
At the top of the discrete variable branch of this $-1$ scheme lie two families of orthogonal
polynomials: the Bannai--Ito (BI) polyno\-mials and their kernel partners the complementary
Bannai--Ito polynomials (CBI); both families correspond to dif\/ferent $q\rightarrow-1$ limits of
the Askey--Wilson polynomials.

The Bannai--Ito polynomials were originally identif\/ied by Bannai and Ito themselves
in~\cite{Ito-1984} where they recognized that these OPs correspond to the $q\rightarrow-1$ limit of
the $q$-Racah polynomials.
However, it is only recently~\cite{Vinet-2012} that the Dunkl shift operator $\mathcal{L}$
admitting the BI polynomials as eigenfunctions has been constructed.
The BI polynomials and their special cases enjoy the Leonard duality property, a~property they
share with all members of the discrete part of the Askey scheme~\cite{Ito-1984,Leonard-1982}.
This means that in addition to satisfying a~three-term recurrence relation, the BI polynomials also
obey a~three-term dif\/ference equation.
From the algebraic point of view, this property corresponds to the existence of an associated
\emph{Leonard pair}~\cite{Terwilliger-2001}.

Amongst the discrete-variable $-1$ polynomials, there are families that do not possess the
Leo\-nard duality property.
That is the case of the complementary Bannai--Ito polynomials and their
descendants~\cite{VZhedanov-2011,Vinet-2012}.
This situation is connected to the fact that in these cases the dif\/ference operator of the
corresponding $q$-polynomials do not admit a~$q\rightarrow-1$ limit.
In~\cite{VZhedanov-2011}, a~\emph{five-term} dif\/ference equation was nevertheless constructed
for the dual $-1$ Hahn polynomials and the def\/ining Dunkl operator for these polynomials was
found.

In this paper, a~one-parameter family of Dunkl operators $\mathcal{D}_{\alpha}$ of which the
complementary Bannai--Ito polynomials are eigenfunctions is derived, thus establishing the
bispectrality of the CBI polynomials.
The operators of this family involve ref\/lections and are of second order in discrete shifts; they
are diagonalized by the CBI polynomials with a~quadratic spectrum.
The corresponding f\/ive-term dif\/ference equation satisf\/ied by the CBI polynomials is presented.
Moreover, an algebra associated to the CBI polynomials is derived.
This quadratic algebra, called the \emph{complementary Bannai--Ito algebra}, is def\/ined in terms
of four generators.
It can be seen as a~deformation with an involution of the quadratic Hahn algebra
${\rm QH}(3)$~\cite{Zhedanov-1992,Zhedanov-1992-2}, which is a~special case of the Askey--Wilson ${\rm AW}(3)$
algebra~\cite{Terwilliger-2011,Zhedanov-1991}.

The paper, which provides a~comprehensive description of the CBI polynomials and their properties,
is organized in the following way.
In Section~\ref{section2}, we present a~review of the Bannai--Ito polynomials.
In Section~\ref{section3}, we def\/ine the complementary Bannai--Ito polynomials and obtain their recurrence
and orthogonality relations.
In Section~\ref{section4}, we use a~proper $q\rightarrow-1$ limit of the Askey--Wilson dif\/ference operator to
construct an operator $\mathcal{D}$ of which the CBI polynomials are eigenfunctions.
We use a~``hidden'' eigenvalue equation to show that one has in fact a~one-parameter family of
operators $\mathcal{D}_{\alpha}$, parametrized by a~complex number $\alpha$, that is diagonalized
by the CBI polynomials.
In Section~\ref{section5}, we derive the CBI algebra and present some aspects of its irreducible representations.
In Section~\ref{section6}, we discuss the relation between the CBI polynomials and three other families of OPs:
the dual~$-1$ Hahn, the para-Krawtchouk and the classical Hahn polynomials; these OP families are
respectively a~limit and two special cases of the CBI polynomials.
We conclude with a~perspective on the continuum limit and an outlook.

\section[Bannai-Ito polynomials]{Bannai--Ito polynomials}\label{section2}

The Bannai--Ito polynomials were introduced in 1984~\cite{Ito-1984} in the complete
classif\/ication of orthogonal polynomials possessing the Leonard duality property (see Section~\ref{section4}).
It was shown that they can be obtained as a~$q\rightarrow-1$ limit of the $q$-Racah polynomials
and some of their properties were derived.
Recently~\cite{Vinet-2012}, it was observed that the BI polynomials also occur as eigensolutions of
a~f\/irst order Dunkl shift operator.
In the following, we review some of the properties of the BI polynomials; we use the presentation
of~\cite{Vinet-2012}.

The monic BI polynomials $B_{n}(x;\rho_1,\rho_2,r_1,r_2)$, denoted $B_{n}(x)$ for notational
convenience, satisfy the three-term recurrence relation
\begin{gather}
\label{Recurrence-Relation}
B_{n+1}(x)+(\rho_1-A_n-C_n)B_{n}(x)+A_{n-1}C_{n}B_{n-1}(x)=xB_n(x),
\end{gather}
with the initial conditions $B_{-1}(x)=0$ and $B_{0}(x)=1$.
The recurrence coef\/f\/icients $A_n$ and $C_n$ are given by
\begin{gather}\label{Coeff-An}
A_n=
\begin{cases}
\dfrac{(n+2\rho_1-2r_1+1)(n+2\rho_1-2r_2+1)}{4(n+g+1)}, & n~\text{even},
\vspace{1mm}\\
\dfrac{(n+2g+1)(n+2\rho_1+2\rho_2+1)}{4(n+g+1)}, & n~\text{odd},
\end{cases}
\\
C_n=
\begin{cases}
-\dfrac{n(n-2r_1-2r_2)}{4(n+g)}, & n~\text{even},
\vspace{1mm}\\
-\dfrac{(n+2\rho_2-2r_2)(n+2\rho_2-2r_1)}{4(n+g)}, & n~\text{odd},
\end{cases}
\end{gather}
where
\begin{gather*}
g=\rho_1+\rho_2-r_1-r_2.
\end{gather*}
It is seen from the above formulas that the positivity condition $u_n=A_{n-1}C_n>0$ cannot be
satisf\/ied for all $n\in\mathbb{N}$~\cite{Chihara-1978}.
Hence it follows that the Bannai--Ito polynomials can only form a~\emph{finite} set of
positive-def\/inite orthogonal polynomials $B_0(x),\ldots,B_N(x)$, which occurs when the ``local''
positivity condition $u_i>0$ for $i\in\{1,\ldots,N\}$ and the truncation conditions $u_0=0$,
$u_{N+1}=0$ are satisf\/ied.
If these conditions are fulf\/illed, the BI polynomials $B_{n}(x)$ satisfy the discrete
orthogonality relation
\begin{gather*}
\sum_{k=0}^{N}w_{k}B_{n}(x_k)B_{m}(x_k)=h_n\delta_{nm},
\end{gather*}
with respect to the positive weight~$w_k$.
The spectral points~$x_k$ are the simple roots of the polynomial $B_{N+1}(x)$.
The explicit formulae for the weight function $w_{k}$ and the grid points~$x_k$ depend on the
realization of the truncation condition~$u_{N+1}=0$.

If $N$ is even, it follows from~\eqref{Coeff-An} that the condition $u_{N+1}=0$ is tantamount to
one of the following requirements:
\begin{alignat*}{3}
&1)~r_1-\rho_1=\frac{N+1}{2},
\qquad &&
2)~r_2-\rho_1=\frac{N+1}{2}, &
\\
&3)~r_1-\rho_2=\frac{N+1}{2},
\qquad &&
4)~r_2-\rho_2=\frac{N+1}{2}. &
\end{alignat*}
For the cases $1)$ and $2)$, the grid points have the expression
\begin{gather}
\label{BI-Grid}
x_k=(-1)^{k}(k/2+\rho_1+1/4)-1/4,
\end{gather}
and the weights take the form
\begin{gather}
\label{Ortho-Weight}
w_{k}=\frac{(-1)^{\nu}}{\ell!}
\frac{(\rho_1-r_1+1/2)_{\ell+\nu}(\rho_1-r_2+1/2)_{\ell+\nu}(\rho_1+\rho_2+1)_{\ell}(2\rho_1+1)_{\ell}}
{(\rho_1+r_1+1/2)_{\ell+\nu}(\rho_1+r_2+1/2)_{\ell+\nu}(\rho_1-\rho_2+1)_{\ell}},
\end{gather}
where one has $k=2\ell+\nu$ with $\nu\in\{0,1\}$ and where $(a)_{n}=a(a+1)\cdots(a+n-1)$ is the
Pochhammer symbol.
For the cases $3)$ and $4)$, the formulae~\eqref{BI-Grid} and~\eqref{Ortho-Weight} hold under the
substitution $\rho_1\leftrightarrow\rho_2$.

If $N$ is odd, it follows from~\eqref{Coeff-An} that the condition $u_{N+1}=0$ is equivalent to
one of the following restrictions:
\begin{gather*}
i)~\rho_1+\rho_2=-\frac{N+1}{2},
\qquad
ii)~r_1+r_2=\frac{N+1}{2},
\qquad
iii)~\rho_1+\rho_2-r_1-r_2=-\frac{N+1}{2}.
\end{gather*}
The condition $iii)$ leads to a~singularity in $u_n$ when $n=(N+1)/2$ and hence only the condi\-tions~$i)$ and~$ii)$ are admissible.
For the case~$i)$, the formulae~\eqref{BI-Grid} and \eqref{Ortho-Weight} hold under the
substitution $\rho_1\leftrightarrow\rho_2$.
For the case $ii)$, the spectral points are given by
\begin{gather*}
x_{k}=(-1)^{k}(r_1-k/2-1/4)-1/4,
\end{gather*}
and the weight function is given by~\eqref{Ortho-Weight} with the substitutions
$(\rho_1,\rho_2,r_1,r_2)\rightarrow-(r_1,r_2,\rho_1,\rho_2)$.

The Bannai--Ito polynomials can be obtained from a~$q\rightarrow-1$ limit of the Askey--Wilson
polynomials and also have the Bannai--Ito algebra as their characteristic algebra
(see~\cite{Genest-2012-1} and~\cite{Vinet-2012}).

\section{CBI polynomials}\label{section3}

In this section we def\/ine the \emph{complementary Bannai--Ito} polynomials through
a~Christof\/fel transformation of the Bannai--Ito polynomials.
We derive their recurrence relation, hypergeometric representation and orthogonality relations from
their kernel properties.

The complementary Bannai--Ito polynomials $I_{n}(x;\rho_1,\rho_2,r_1,r_2)$, denoted $I_{n}(x)$ for
convenience, are def\/ined from the BI polynomials $B_n(x)$ by the transformation~\cite{Vinet-2012}
\begin{gather}
\label{CBI-Def}
I_{n}(x)=\frac{B_{n+1}(x)-A_n B_{n}(x)}{x-\rho_1},
\end{gather}
where $A_n$ is as in~\eqref{Coeff-An}.
The transformation~\eqref{CBI-Def} is an example of a~Christof\/fel
transformation~\cite{Szego-1975}.
It is easily seen from the def\/inition~\eqref{CBI-Def} that $I_n(x)$ is a~monic polynomial of
degree $n$ in $x$.
The inverse relation for the CBI polynomials is given by a~Geronimus~\cite{Zhedanov-1997}
transformation and has the expression
\begin{gather}
\label{CBI-Inv}
B_{n}(x)=I_{n}(x)-C_nI_{n-1}(x).
\end{gather}
This formula can be verif\/ied by direct substitution of~\eqref{CBI-Def} in~\eqref{CBI-Inv} which
yields back the def\/ining relation~\eqref{Recurrence-Relation} of the BI polynomials.
In the reverse, the substitution of~\eqref{CBI-Inv} in~\eqref{CBI-Def} yields the three-term
recurrence relation~\cite{Karlin-1957}
\begin{gather}
\label{Recurrence-1}
I_{n+1}(x)+(\rho_1-A_n-C_{n+1})I_{n}(x)+A_{n}C_{n} I_{n-1}(x)=xI_{n}(x),
\end{gather}
where $A_n$ and $C_n$ are given by~\eqref{Coeff-An}.
The recurrence relation~\eqref{Recurrence-1} can be written explicitly as
\begin{gather}
\label{CBI-Off}
I_{n+1}(x)+(-1)^{n}\rho_2 I_{n}(x)+\tau_{n}I_{n-1}(x)=xI_{n}(x),
\end{gather}
where $\tau_n$ is given by
\begin{gather}
\tau_{2n}=-\frac{n(n+\rho_1-r_1+1/2)(n+\rho_1-r_2+1/2)(n-r_1-r_2)}{(2n+g)(2n+g+1)},
\nonumber
\\
\tau_{2n+1}=-\frac{(n+g+1)(n+\rho_1+\rho_2+1)(n+\rho_2-r_1+1/2)(n+\rho_2-r_2+1/2)}{(2n+g+1)(2n+g+2)},
\label{Recurrence--CBI}
\end{gather}
and where $g=\rho_1+\rho_2-r_1-r_2$.
One has also the initial conditions $I_0=1$ and $I_1=x-\rho_2$.
The CBI polynomials are kernel polynomials of the BI polynomials.
Indeed, by noting that
\begin{gather*}
A_n=B_{n+1}(\rho_1)/B_n(\rho_1),
\end{gather*}
which follows by induction from~\eqref{Recurrence-Relation}, the transformation~\eqref{CBI-Def} may
be cast in the form
\begin{gather}
\label{Kernel}
I_{n}(x)=(x-\rho_1)^{-1}\left[B_{n+1}(x)-\frac{B_{n+1}(\rho_1)}{B_n(\rho_1)}B_n(x)\right].
\end{gather}
It is manifest from~\eqref{Kernel} that $I_n(x)$ are the kernel polynomials associated to $B_n(x)$
with kernel parameter $\rho_1$~\cite{Chihara-1978}.
Since the BI polynomials $B_{n}(x)$ are orthogonal with respect to a~linear functional
$\sigma^{(i)}$:
\begin{gather*}
\langle\sigma^{(i)},B_n(x)B_m(x)\rangle=0,
\qquad
n\neq m,
\end{gather*}
where the upper index on $\sigma^{(i)}$ designates the possible functionals associated to the
various truncation conditions, it follows from~\eqref{Kernel} that we have~\cite{Chihara-1978}
\begin{gather}
\label{Kernel-Ortho}
\langle\sigma^{(i)},(x-\rho_1)I_n(x)I_m(x)\rangle=0,
\qquad
n\neq m.
\end{gather}
Hence the orthogonality and positive-def\/initeness of the CBI polynomials can be studied using the
formulae~\eqref{Recurrence--CBI} and \eqref{Kernel-Ortho}.

It is seen from~\eqref{Recurrence--CBI} that the condition $\tau_n>0$ cannot be ensured for all $n$
and hence the complementary Bannai--Ito polynomials can only form a~\emph{finite} system of
positive-def\/inite orthogonal polynomials $I_0(x),\ldots,I_N(x)$, provided that the ``local''
positivity $\tau_n>0$, $n\in\{1,\ldots,N\}$, and truncation conditions $\tau_0=0$
and $\tau_{N+1}=0$ are satisf\/ied.

When $N$ is even, the truncation conditions $\tau_0=0$ and $\tau_{N+1}=0$ are equivalent to one of
the four prescriptions
\begin{alignat}{3}
&1)~\rho_2-r_1=-\frac{N+1}{2},
\qquad &&
2)~\rho_2-r_2=-\frac{N+1}{2}, &
\nonumber
\\
&3)~\rho_1+\rho_2=-\frac{N+2}{2},
\qquad &&
4)~g=-\frac{N+2}{2}. &
\label{CBI-Trunca-Even}
\end{alignat}
Since the condition $4)$ leads to a~singularity in $\tau_n$, only the conditions $1)$, $2)$
and $3)$ are admis\-sible.
For all three conditions and assuming that the positivity conditions are satisf\/ied, the CBI
polynomials enjoy the orthogonality relation
\begin{gather}
\label{Ortho-2}
\sum_{k=0}^{N}\widetilde{w}_{k}I_{n}(x_k)I_m(x_k)=\widetilde{h}_n\delta_{nm},
\end{gather}
where the spectral points are given by
\begin{gather*}
x_{k}=(-1)^{k}(k/2+\rho_2+1/4)-1/4,
\end{gather*}
and the positive weights are
\begin{gather*}
\widetilde{w}_{k}=(x_k-\rho_1)w_k,
\end{gather*}
with $w_k$ def\/ined by~\eqref{Ortho-Weight} with the substitution $\rho_1\leftrightarrow\rho_2$.

When $N$ is odd, the truncation conditions $\tau_0=0$ and $\tau_{N+1}=0$ are tantamount to
\begin{gather}
\label{CBI-Trunca-Odd}
i)~r_1-\rho_1=\frac{N+2}{2},
\qquad
ii)~r_1+r_2=\frac{N+1}{2},
\qquad
iii)~r_2-\rho_1=\frac{N+2}{2}.
\end{gather}
If the positivity condition $\tau_n>0$ is satisf\/ied for $n\in\{1,\ldots,N\}$, the CBI polynomials
will enjoy the orthogonality relation~\eqref{Ortho-2} with respect to the positive def\/inite
weight function $\widetilde{w}_k$.
When either condition $i)$ or $ii)$ is satisf\/ied, the spectral points are given by
\begin{gather*}
x_k=(-1)^{k}(r_1-k/2-1/4)-1/4,
\end{gather*}
together with the weight function $\widetilde{w}_k=(x_k-\rho_1)w_k$ where $w_k$ is given
by~\eqref{Ortho-Weight} with the replacement $(\rho_1,\rho_2,r_1,r_2)=-(r_1,r_2,\rho_1,\rho_2)$.
Finally, the orthogonality relation for the truncation condition $iii)$ is obtained from the
preceding case under the exchange $r_1\leftrightarrow r_2$.

Let us now illustrate when positive-def\/initeness occurs for the CBI polynomials.
We f\/irst consider the even $N$ case.
It is suf\/f\/icient to take
\begin{gather}
\label{Para-Even}
\rho_1=\left(\frac{\frac{a+b}{2}+c+N}{2}\right)\!,
\quad
\rho_2=\left(\frac{\frac{a+b}{2}-1}{2}\right)\!,
\quad
r_1=\left(\frac{\frac{a+b}{2}+N}{2}\right)\!,
\quad
r_2=\left(\frac{a-b}{4}\right)\!,
\end{gather}
where $a$, $b$ and $c$ are arbitrary positive parameters.
Assuming~\eqref{Para-Even}, the recurrence coef\/f\/i\-cients~\eqref{Recurrence--CBI} become
\begin{gather}\label{voila}
\tau_n=
\begin{cases}
\dfrac{n(N-n+a)(n+c+1)(n+b+c+N+1)}{16(n+g)(n+g+1)},&n~\text{even},
\vspace{1mm}\\
\dfrac{(N-n+1)(n+b-1)(n+b+c)(n+a+b+c+N)}{16(n+g)(n+g+1)},&n~\text{odd},
\end{cases}
\end{gather}
where $g=(b+c-1)/2$.
It is obvious from~\eqref{voila} that the positivity and truncation conditions are satisf\/ied for
$n\in\{1,\ldots,N\}$; this corresponds to the case $1)$ of~\eqref{CBI-Trunca-Even}.

Consider the situation when $N>1$ is odd.
We introduce the parametrization
\begin{gather}\label{Para-Odd}
\rho_1\!=\!\left(\frac{\frac{\zeta+\xi}{2}+\chi+N}{2}\right)\!,
\quad
\rho_2\!=\!\left(\frac{\zeta-\xi}{4}\right)\!,
\quad
r_1\!=\!\left(\frac{\frac{\zeta+\xi}{2}+N+1}{2}\right)\!,
\quad
r_2\!=\!-\!\left(\frac{\zeta+\xi}{4}\right)\!,
\end{gather}
where $\zeta$, $\xi$ and $\chi$ are arbitrary positive parameters.
The recurrence coef\/f\/icients become
\begin{gather*}
\tau_{n}=
\begin{cases}
\dfrac{n(N-n+1)(n+\chi)(n+\zeta+\xi+\chi+N+1)}{16(n+g)(n+g+1)}, & n~\text{even},
\vspace{1mm}\\
\dfrac{(N-n+\xi+1)(n+\zeta)(n+\zeta+\chi)(n+\zeta+\chi+N+1)}{16(n+g)(n+g+1)}, & n~\text{odd},
\end{cases}
\end{gather*}
with $g=(\zeta+\chi-1)/2$.
Assuming~\eqref{Para-Odd}, the positivity and truncation conditions are manifestly fulf\/illed;
this corresponds to the condition $ii)$ of~\eqref{CBI-Trunca-Odd}.
The other cases can be treated in similar fashion.

It is possible to derive a~hypergeometric representation for the CBI polynomials using
a~me\-thod~\cite{Vinet-2012, Vinet-2011-2} which is analogous to Chihara's construction of symmetric
orthogonal polyno\-mials~\cite{Chihara-1978} and closely related to the scheme developed
in~\cite{Marcellan-1997} (see also~\cite{Chihara-1964}).
Given the three-term recurrence relation~\eqref{CBI-Off}, it follows by induction that the
polynomials~$I_{n}(x)$ can be written as
\begin{gather}
\label{vlan}
I_{2n}=R_{n}\big(x^2\big),
\qquad
I_{2n+1}=(x-\rho_2)Q_{n}\big(x^2\big),
\end{gather}
where $R_{n}(x^2)$ and $Q_{n}(x^2)$ are monic polynomials of degree~$n$.
It follows directly from~\eqref{vlan} and~\eqref{CBI-Off} that the polynomials~$R_{n}(x^2)$
and~$Q_{n}(x^2)$ obey the following system of recurrence relations
\begin{gather*}
R_{n}(z)=Q_n(z)+\tau_{2n}Q_{n-1}(z),
\qquad
\big(z-\rho_2^2\big)Q_n(z)=R_{n+1}(z)+\tau_{2n+1}R_{n}(z).
\end{gather*}
This system is equivalent to the following pair of equations:
\begin{gather*}
R_{n+1}(z)+\big(\rho_2^2+\tau_{2n}+\tau_{2n+1}\big)R_n(z)+\tau_{2n-1}\tau_{2n}R_{n-1}(z)=zR_{n}(z),
\\
Q_{n+1}(z)+\big(\rho_2^2+\tau_{2n+1}+\tau_{2n+2}\big)Q_{n}(z)+\tau_{2n}\tau_{2n+1}Q_{n-1}(z)= zQ_{n}(z).
\end{gather*}
These recurrence relations can be identif\/ied with those of the Wilson
polynomials~\cite{Koekoek-2010}.
From this identif\/ication, we obtain
\begin{gather}
R_{n}\big(x^2\big)=\eta_n\;
\pFq{4}{3}{-n,n+g+1,\rho_2+x,\rho_2-x}{\rho_1+\rho_2+1,\rho_2-r_1+1/2,\rho_2-r_2+1/2}{1},
\nonumber
\\
Q_{n}\big(x^2\big)=
\iota_{n}\;\pFq{4}{3}{-n,n+g+2,\rho_2+1+x,\rho_2+1-x}{\rho_1+\rho_2+2,\rho_2-r_1+3/2,\rho_2-r_2+3/2}{1},
\label{CBI-Hypergeo}
\end{gather}
where ${}_pF_{q}$ denotes the generalized hypergeometric function~\cite{Gasper-1990} and where the
normalization coef\/f\/icients, which ensure that the polynomials are monic, are given by
\begin{gather*}
\eta_n=\frac{(\rho_1+\rho_2+1)_n(\rho_2-r_1+1/2)_n(\rho_2-r_2+1/2)_n}{(n+g+1)_{n}},
\\
\iota_n =\frac{(\rho_1+\rho_2+2)_n(\rho_2-r_1+3/2)_n(\rho_2-r_2+3/2)_n}{(n+g+2)_{n}}.
\end{gather*}
Thus the monic CBI polynomials have the hypergeometric representation~\eqref{vlan}.
For def\/initeness and future reference, let us now gather the preceding results in the following
proposition.
\begin{proposition}
The complementary Bannai--Ito polynomials $I_{n}(x;\rho_1,\rho_2,r_1,r_2)$ are the kernel
polynomials of the Bannai--Ito polynomials $B_{n}(x;\rho_1,\rho_2,r_1,r_2)$ with kernel parameter
$\rho_1$.
The monic CBI polynomials obey the three-term recurrence relation
\begin{gather*}
I_{n+1}(x)+(-1)^{n}\rho_2I_{n}(x)+\tau_{n}I_{n-1}(x)=xI_{n}(x),
\end{gather*}
where $\tau_{n}$ is given by~\eqref{Recurrence--CBI}.
They have the explicit hypergeometric representation
\begin{gather*}
I_{2n}(x)=R_{n}\big(x^2\big),
\qquad
I_{2n+1}(x)=(x-\rho_2)Q_{n}\big(x^2\big),
\end{gather*}
where $R_{n}(x^2)$ and $Q_{n}(x^2)$ are as specified by~\eqref{CBI-Hypergeo}.
If the truncation condition $\tau_{N+1}=0$ and the positivity condition $\tau_{n}>0$,
$n\in\{1,\ldots,N\}$, are satisfied, the CBI polynomials obey the orthogonality relation
\begin{gather*}
\sum_{k=0}^{N}\widetilde{w}_{k}I_{n}(x_k)I_{m}(x_k)=\widetilde{h}_{n}\delta_{nm},
\end{gather*}
with respect to the positive weights $\widetilde{w}_k$.
The grid points $x_k$ correspond to the simple roots of the polynomial $I_{N+1}(x)$.
The formulas for the weights and grid points depend on the truncation condition.
With $w_{k}(\rho_1,\rho_2,r_1,r_2)$ given as in~\eqref{Ortho-Weight}, one has
\begin{enumerate}\itemsep=0pt
\item For $r_1=\frac{N+1}{2}+\rho_2$, $r_2=\frac{N+1}{2}+\rho_2$ or $\rho_1=-\frac{N+2}{2}-\rho_2$
with $N$ even:
\begin{gather*}
x_k=(-1)^{k}(\rho_2+k/2+1/4)-1/4,
\qquad
\widetilde{w}_{k}=(x_k-\rho_1)w_k(\rho_2,\rho_1,r_1,r_2).
\end{gather*}
\item For $r_1=\frac{N+2}{2}+\rho_1$ or $r_1=\frac{N+1}{2}-r_2$ with $N$ odd:
\begin{gather*}
x_{k}=(-1)^{k}(r_1-k/2-1/4)-1/4,
\qquad
\widetilde{w}_{k}=(x_k-\rho_1)w_{k}(-r_1,-r_2,-\rho_1,-\rho_2).
\end{gather*}
\item For $r_2=\frac{N+2}{2}+\rho_1$ with $N$ odd:
\begin{gather*}
x_{k}=(-1)^{k}(r_2-k/2-1/4)-1/4,
\qquad
\widetilde{w}_{k}=(x_k-\rho_1)w_{k}(-r_2,-r_1,-\rho_1,-\rho_2).
\end{gather*}
\end{enumerate}
\end{proposition}

\begin{proof}
The proof follows from the above considerations.
\end{proof}

Note that the normalization factor $\widetilde{h}_{n}$ appearing in~\eqref{Ortho-2} can easily be
evaluated in terms of the product $\tau_1\tau_2\cdots\tau_n$.

The complementary Bannai--Ito polynomials can be obtained from the Askey--Wilson polynomials upon
taking the $q\rightarrow-1$ limit~\cite{Vinet-2012}.
Consider the Askey--Wilson polynomials~\cite{Koekoek-2010} $p_{n}(z;a,b,c,d)$
\begin{gather}
\label{Def-AW}
p_n(z;a,b,c,d)=a^{-n}(ab,ac,ad;q)_{n}\;\pfq{4}{3}{q^{-n},abcdq^{n-1},az,az^{-1}}{ab,ac,ad}{q;q},
\end{gather}
where $\phi$ denotes the basic generalized hypergeometric function~\cite{Gasper-1990}.
These polynomials depend on the argument $x=(z+z^{-1})/2$ and on four complex parameters $a$, $b$,
$c$ and $d$.
They obey the recurrence relation~\cite{Koekoek-2010}
\begin{gather}
\label{Recu-AW}
\alpha_{n}p_{n+1}(z)+\big(a+a^{-1}-\alpha_n-\gamma_n\big)p_{n}(z)+\gamma_n p_{n-1}(z)=
\big(z+z^{-1}\big)p_n(z),
\end{gather}
where the coef\/f\/icients are
\begin{gather*}
\alpha_n=\frac{(1-abq^{n})(1-acq^{n})(1-adq^{n})(1-abcdq^{n-1})}{a(1-abcdq^{2n-1})(1-abcdq^{2n})},
\\
\gamma_n=\frac{a(1-q^{n})(1-bcq^{n-1})(1-bdq^{n-1})(1-cdq^{n-1})}{(1-abcdq^{2n-2})(1-abcdq^{2n-1})}.
\end{gather*}
To recover the CBI polynomials, we consider the parametrization
\begin{alignat}{4}
& a=ie^{\epsilon (2\rho_1+3/2)},
\qquad &&
b=-ie^{\epsilon (2\rho_2+1/2)},
\qquad &&
c=ie^{\epsilon (-2r_2+1/2)},&
\nonumber
\\
& d=ie^{\epsilon (-2r_1+1/2)},
\qquad &&
q=-e^{\epsilon},
\qquad &&
z=ie^{-2\epsilon y}.&
\label{Para-Limit}
\end{alignat}
It can be verif\/ied that the limit $q\rightarrow-1$ of the Askey--Wilson polynomials
\begin{gather*}
\lim_{q\rightarrow-1}p_{n}(z)=p_n^*(y),
\end{gather*}
exists~\cite{Vinet-2012} and that $p_n^*(y)$ is a~polynomial of degree~$n$ in the variable~$y$.
Dividing the recurrence relation~\eqref{Recu-AW} by $1+q$ and taking the limit
$\epsilon\rightarrow0$, which amounts to taking $q\rightarrow-1$, one f\/inds that the recurrence
relation of the limit polynomials $p_{n}^*(y)$ is
\begin{gather*}
\alpha_n^{*}p_{n+1}^{*}(y)+(-1)^{n}\rho_2 p_{n}^{*}(y)+\gamma_n^{*}p_{n-1}^{*}(y)=
(y-1/4)p_{n}^{*}(y),
\end{gather*}
where
\begin{gather}
\alpha_{2n}^{*}=-\frac{(n+\rho_1+\rho_2+1)(n+g+1)}{(2n+g+1)},
\nonumber
\\
\alpha_{2n+1}^{*}=
-\frac{(n+\rho_1-r_1+3/2)(n+\rho_1-r_2+3/2)}{(2n+g+2)},
\nonumber
\\
\gamma_{2n}^{*}=\frac{n(n-r_1-r_2)}{(2n+g+1)},
\nonumber
\\
\gamma_{2n+1}^{*}=
\frac{(n+\rho_2-r_1+1/2)(n+\rho_2-r_2+1/2)}{(2n+g+2)}.
\label{Recu-Lim}
\end{gather}
From~\eqref{Recurrence--CBI} and \eqref{Recu-Lim}, one has the identif\/ication
\begin{gather}
\label{CBI-Identification}
\widehat{p_{n}^{*}}(y)=I_n(y-1/4),
\end{gather}
where $\widehat{p_{n}^{*}}$ are the monic version of the limit polynomials $p_{n}^{*}(y)$.
Consequently, the CBI polynomial correspond to a~$q\rightarrow-1$ limit of the Askey--Wilson
polynomials, up to a~shift in argument.
This property will be used in the next section to construct a~Dunkl operator that has the CBI
polynomials as eigenfunctions.

\section{Bispectrality of CBI polynomials}\label{section4}

In this section, we obtain a~family of second order Dunkl shift operators for which the
complementary Bannai--Ito polynomials are eigenfunctions with eigenvalues quadratic in $n$.
This family will be constructed from a~limit of a~quadratic combination of the Askey--Wilson
$q$-dif\/ference operator.
We shall refer to these operators as the def\/ining operators of the CBI polynomials.

Consider the Askey--Wilson polynomials $p_n(x)$ def\/ined by~\eqref{Def-AW}.
They obey the $q$-dif\/ference equation~\cite{Koekoek-2010}
\begin{gather}
\label{AW-Operator}
\big(\Omega(z)E_{z,q}+\Omega\big(z^{-1}\big)E_{z,q^{-1}}-\big(\Omega(z)+\Omega\big(z^{-1}\big)\big)\mathbb{I}\big)p_{n}(z)
=\Lambda_{n}p_{n}(z),
\end{gather}
where $E_{z,q}f(z)=f(qz)$ is the $q$-shift operator and $\mathbb{I}$ denotes the identity.
The l.h.s.\ of~\eqref{AW-Operator} is the Askey--Wilson operator.
The eigenvalues take the form
\begin{gather*}
\Lambda_n=\big(q^{-n}-1\big)\big(1-abcd q^{n-1}\big),
\end{gather*}
and the coef\/f\/icient $\Omega(z)$ is given by
\begin{gather*}
\Omega(z)=\frac{(1-az)(1-bz)(1-cz)(1-dz)}{(1-z^2)(1-qz^2)}.
\end{gather*}
We now consider the limiting form of the $q$-dif\/ference equation~\eqref{AW-Operator} when
$q\rightarrow-1$.
As done previously, we choose the parametrization~\eqref{Para-Limit}, which correspond to the CBI
polynomials.
We already showed that the Askey--Wilson polynomials $p_n(z)$ become the complementary Bannai--Ito
polynomials $p_n^{*}(y)$.
In the limit $q\rightarrow-1$, the $q$-shift operation $p_{n}(z)\rightarrow p_{n}(qz)$ becomes
$p_n^{*}(y)\rightarrow p_{n}^{*}(-y+1/2)$ while $p_n(z)\rightarrow p_n(q^{-1}z)$ is reduced to
$p_n^{*}(y)\rightarrow p_n^{*}(-y-1/2)$.

It is natural to expect that in the limit $q\rightarrow-1$, the equation~\eqref{AW-Operator} will
yield a~def\/ining operator for the CBI polynomials.
However a~direct computation shows that the limit $\epsilon\rightarrow0$ of the
equation~\eqref{AW-Operator} with the parametrization~\eqref{Para-Limit} \emph{does not exist}.
It is hence impossible to f\/ind the desired operator for the CBI polynomials directly from
a~limiting procedure on equation~\eqref{AW-Operator}.
Nevertheless, it is possible to work around this dif\/f\/iculty by choosing an appropriate
\emph{quadratic} combination of the Askey--Wilson operator that survives the limit $q\rightarrow-1$.
A similar procedure was used in~\cite{VZhedanov-2011} to establish the bispectrality of the dual
$-1$ Hahn polynomials.

Let $\mathcal{O}$ denote the Askey--Wilson operator
\begin{gather*}
\mathcal{O}=
\Omega(z)E_{z,q}+\Omega\big(z^{-1}\big)E_{z,q^{-1}}-\big(\Omega(z)+\Omega\big(z^{-1}\big)\big)\mathbb{I},
\end{gather*}
which acts on the space of functions $f(z)$ of argument $z$.
We consider the following quadratic combination
\begin{gather*}
\mathcal{T}=c^{(2)}\mathcal{O}^2+c^{(1)}\mathcal{O},
\end{gather*}
with
\begin{gather*}
c^{(2)}=\frac{1}{16(1+q)^2},
\qquad
c^{(1)}=\frac{1}{4}\left(\frac{1}{(q+1)^2}-\frac{g+1}{q+1}\right),
\end{gather*}
where $g=\rho_1+\rho_2-r_1-r_2$.
Since the operator $\mathcal{O}$ acts diagonally on the Askey--Wilson polynomials, we have
\begin{gather}\label{Startup}
\mathcal{T}p_{n}(z)=\big(c^{(2)}\Lambda_{n}^2+c^{(1)}\Lambda_{n}\big)p_n(x).
\end{gather}
Upon taking the limit $\epsilon\rightarrow0$ with the parametrization~\eqref{Para-Limit}, the
relation~\eqref{Startup} becomes{\samepage
\begin{gather}
\Phi_1(y)p_{n}^{*}(y+1)+\big\{\Phi_5(y)-\Phi_2(y)-\Phi_3(y)\big\}p_{n}^{*}(1/2-y)
+\big\{\Phi_3(y)-\Phi_4(y)-\Phi_5(y)\big\}p_n^{*}(y)
\nonumber
\\
\qquad
{}+\big\{\Phi_4(y)-\Phi_1(y)\big\}p_{n}^{*}(-y-1/2)+\Phi_2(y)p_n^{*}(y-1)=\kappa_{n}p_n^{*}(y)
\label{Dunkl-1}
\end{gather}
where} the eigenvalues are
\begin{gather*}
\kappa_{2n}=n^2+(g+1)n,
\qquad
\kappa_{2n+1}=n^2+(g+2)n+g^2+2g+5/4.
\end{gather*}
The coef\/f\/icients $\Phi_i(y)$ are given by
\begin{gather*}
\Phi_1(y) =\frac{(y+\rho_1+3/4)(y+\rho_2+3/4)(y-r_1+1/4)(y-r_2+1/4)}{4(y+1/4)(y+3/4)},
\\
\Phi_2(y) =\frac{(y-\rho_1-5/4)(y-\rho_2-1/4)(y+r_1-3/4)(y+r_2-3/4)}{4(y-1/4)(y-3/4)},
\\
\Phi_3(y) =\frac{(y+\rho_1+3/4)(y-\rho_2-1/4)(y-r_1+1/4)(y-r_2+1/4)}{4(y-1/4)(y+1/4)},
\\
\Phi_4(y) =\frac{(y+\rho_1+3/4)(y+\rho_2-1/4)(y-r_1+1/4)(y-r_2+1/4)}{4(y-1/4)(y+1/4)},
\\
\Phi_5(y) =\frac{(y-\rho_2-1/4)}{4(y-1/4)}\big\{2y^2-y+\nu\big\},
\end{gather*}
where $\nu$ takes the form
\begin{gather*}
\nu=r_1+r_2+2r_1r_2-2\rho_1-2(r_1+r_2)\rho_1-4\rho_2+1/8-2g^2.
\end{gather*}
By the identif\/ication~\eqref{CBI-Identification}, the relation~\eqref{Dunkl-1} gives the
complementary Bannai--Ito polynomials $I_n(x)$ as eigenfunctions of a~second order Dunkl shift
operator, hence establishing their bispectrality property.
In operator form, the equation~\eqref{Dunkl-1} may be rewritten as
\begin{gather*}
\mathcal{H} I_{n}(y-1/4)=\kappa_n I_{n}(y-1/4),
\end{gather*}
where $\mathcal{H}$ has the expression
\begin{gather*}
\mathcal{H}=\Phi_1T^{1}+(\Phi_4-\Phi_1)T^{1/2}R+(\Phi_3-\Phi_4-\Phi_5)\mathbb{I}
+(\Phi_5-\Phi_2-\Phi_3)T^{-1/2}R+\Phi_2T^{-1},
\end{gather*}
where $T^{h}f(y)=f(y+h)$ and $Rf(y)=f(-y)$.
Upon applying the unitary transformation
\begin{gather*}
\widetilde{\mathcal{H}}=T^{1/4}\mathcal{H}T^{-1/4},
\end{gather*}
on the operator $\mathcal{H}$ and changing the variable from $y$ to $x$, the eigenvalue
equation~\eqref{Dunkl-1} for the CBI polynomials becomes
\begin{gather*}
\widetilde{\mathcal{H}}I_n(x)=\kappa_n I_n(x),
\end{gather*}
where we have
\begin{gather}
\label{eq}
\widetilde{\mathcal{H}}=\widetilde{\Phi}_1
T^{+}+(\widetilde{\Phi}_4-\widetilde{\Phi}_1)T^{+}R
+(\widetilde{\Phi}_3-\widetilde{\Phi}_4-\widetilde{\Phi}_5)\mathbb{I}
+(\widetilde{\Phi}_5-\widetilde{\Phi}_2-\widetilde{\Phi}_3)R+\widetilde{\Phi}_2T^{-},
\end{gather}
with $T^{+}=T^{1}$ and $T^{-}=T^{-}$ the usual shift operators in $x$.
The coef\/f\/icients now have the expression
\begin{gather*}
\widetilde{\Phi}_i=\Phi_i(x+1/4).
\end{gather*}

We now turn to the study of the uniqueness of the operator $\mathcal{H}$ which def\/ines the
eigenvalue equation of the complementary Bannai--Ito polynomials (apart from trivial af\/f\/ine
transformations).
Quite strikingly, a~one-parameter family of such operators can be constructed.
This peculiarity is due to the presence of a~``hidden'' symmetry in the CBI polynomials.
To see this, we recall the relation~\eqref{vlan} for the CBI polynomials
\begin{gather*}
I_{2n}=R_{n}\big(x^2\big),
\qquad
I_{2n+1}=(x-\rho_2)Q_{n}\big(x^2\big),
\end{gather*}
where $R_n(x^2)$ and $Q_{n}(x^2)$ are monic polynomials of degree $n$.
From the above relation, it is easily seen that
\begin{gather*}
I_{2n}(-x)=I_{2n}(x),
\qquad
\text{and}
\qquad
I_{2n+1}(-x)=\frac{(x+\rho_2)}{(\rho_2-x)}I_{2n+1}(x).
\end{gather*}
The above equations are equivalent to the following non-trivial ``hidden" eigenvalue equation for
the CBI polynomials
\begin{gather*}
\frac{(\rho_2-x)}{2x}\big(I_{n}(-x)-I_n(x)\big)=\mu_n I_{n}(x),
\end{gather*}
where $\mu_{2n}=0$ and $\mu_{2n+1}=1$.
In operator form, we write
\begin{gather}
\label{Hidden-Eigen}
\frac{(x-\rho_2)}{2x} (\mathbb{I}-R )I_n(x)=\mathcal{U}I_n(x)=\mu_n I_n(x).
\end{gather}
The equation~\eqref{Hidden-Eigen} indicates that adding $\alpha \mathcal{U}$ to the
operator~\eqref{eq} will give another eigenvalue equation for the complementary Bannai--Ito
polynomials.
The modif\/ied operator
\begin{gather*}
\widetilde{\mathcal{H}}'=\widetilde{\mathcal{H}}+\alpha\mathcal{U},
\end{gather*}
will have the same spectrum as $\widetilde{\mathcal{H}}$ in the even sector; in the odd sector, the
eigenvalues will dif\/fer by the constant parameter $\alpha$.

For def\/initeness and future reference, let us now collect the preceding results in the following
theorem.
\begin{theorem}
Let $\mathcal{D}_0$ be the second order Dunkl shift operator acting on the space of functions~$f(x)$ of argument~$x$
\begin{gather}
\mathcal{D}_0=A(x) T^{+}+B(x) T^{-}+C(x) R+D(x) T^{+}R
-(A(x)+B(x)+C(x)+D(x))\mathbb{I},\label{CBI-D0}
\end{gather}
where $T^{\pm}f(x)=f(x\pm1)$ and $Rf(x)=f(-x)$, with the coefficients
\begin{gather*}
A(x)=\frac{(x+\rho_1+1)(x+\rho_2+1)(2x-2r_1+1)(2x-2r_2+1)}{8(x+1)(2x+1)},
\\
B(x)=\frac{(x-\rho_2)(x-\rho_1-1)(2x+2r_1-1)(2x+2r_2-1)}{8x(2x-1)},
\\
C(x)=\frac{(x-\rho_2)(4x^2+\omega)}{8x}-\frac{(x-\rho_2)(x+\rho_1+1)(2x-2r_1+1)(2x-2r_2+1)}{8x(2x+1)}-B(x),
\\
D(x)=\frac{\rho_2(x+\rho_1+1)(2x-2r_1+1)(2x-2r_2+1)}{8x(x+1)(2x+1)},
\end{gather*}
and with
\begin{gather*}
\omega=4\rho_1-4(r_1+r_2)\rho_1+4r_1r_2-6(r_1+r_2)+5.
\end{gather*}
Furthermore, let $\alpha\in\mathbb{C}$ be a~complex number and denote the monic complementary
Bannai--Ito polynomials by~$I_n(x)$.
Then the following eigenvalue equation is satisfied:
\begin{gather}
\label{Eigen-Equation}
\mathcal{D}_{\alpha}I_{n}(x)=\Lambda_{n}^{(\alpha)}I_n(x),
\end{gather}
where the eigenvalues are
\begin{gather}
\label{Eigenvalues}
\Lambda_{2n}^{(\alpha)}=n^2+(g+1)n,
\qquad
\Lambda_{2n+1}^{(\alpha)}=n^2+(g+2)n+\alpha,
\end{gather}
and where we have defined
\begin{gather}
\label{Full--OP}
\mathcal{D}_{\alpha}=\mathcal{D}_0+\alpha \frac{(x-\rho_2)}{2x}(\mathbb{I}-R).
\end{gather}
\end{theorem}

\begin{proof}
The result follows from the above considerations.
\end{proof}

We now discuss the CBI polynomials in the context of the Leonard duality.
A family of orthogonal polynomials $\mathcal{P}_n(x)$ is said to possess the \emph{Leonard duality
property} if it satisf\/ies both a~three-term recurrence relation with respect to $n$
and a~three-term dif\/ference equation of the form
\begin{gather*}
\theta(x_k)\mathcal{P}_n(x_{k+1})+\nu(x_k)\mathcal{P}_n(x_k)+\mu(x_k)\mathcal{P}_n(x_{k-1})=
\vartheta_{n}\mathcal{P}_n(x_k),
\end{gather*}
on a~discrete set of points $x_k$, $k\in\mathbb{Z}$.
The classif\/ication of the polynomials with this property was f\/irst accomplished by Leonard
in~\cite{Leonard-1982}; his theorem was later generalized to include inf\/inite-dimensional grids
by Bannai and Ito~\cite{Ito-1984}.
It turns out that the complementary Bannai--Ito polynomials lie beyond the scope of the Leonard
duality.
Indeed, the operators $\mathcal{D}_{\alpha}$ can be used to show that the CBI polynomials obey
a~\emph{five-term} dif\/ference equation on an inf\/inite-dimensional grid.
This result is obtained in the following way.

First consider the grid $x_k$ def\/ined by
\begin{gather}
\label{BI-Grid-1}
x_k=(-1)^{k}(k/2+h+1/4)-1/4,
\qquad
k\in\mathbb{Z},
\end{gather}
where $h$ is an arbitrary real parameter.
It is easily seen that the grid~\eqref{BI-Grid-1} is preserved by the operators appearing
in~\eqref{CBI-D0}.
Explicitly, we have
\begin{alignat*}{3}
& T^{+}x_{k} =
\begin{cases}
x_{k+2}, & k~\text{even},
\\
x_{k-2}, & k~\text{odd},
\end{cases}
\qquad &&
T^{-}x_{k}=
\begin{cases}
x_{k-2}, & k~\text{even},
\\
x_{k+2}, & k~\text{odd},
\end{cases} &
\\
& Rx_{k} =
\begin{cases}
x_{k-1}, & k~\text{even},
\\
x_{k+1}, & k~\text{odd},
\end{cases}
\qquad &&
T^{+}R x_k=
\begin{cases}
x_{k+1}, & k~\text{even},
\\
x_{k-1}, & k~\text{odd}.
\end{cases}&
\end{alignat*}
Referring to $\mathcal{D}_0$, one f\/inds the following f\/ive-term dif\/ference equation for the
CBI polynomials:
\begin{gather}
\nonumber
u(x_k)I_n(x_{k+2})+v(x_k)I_n(x_{k+1})+m(x_k)I_n(x_{k})
\\
\qquad
{} +t(x_k)I_n(x_{k-1})+r(x_k)I_{m}(x_{k-2})=\Lambda_{n}^{(0)} I_{n}(x_k)
\label{CBI-Difference}
\end{gather}
where we have
\begin{alignat*}{3}
& u(x_{k}) =
\begin{cases}
A(x_k), & k~\text{even},
\\
B(x_k), & k~\text{odd},
\end{cases}
\qquad &&
v(x_k)=
\begin{cases}
D(x_k), & k~\text{even},
\\
C(x_k), & k~\text{odd},
\end{cases} &
\\
& t(x_k) =
\begin{cases}
C(x_k), & k~\text{even},
\\
D(x_k), & k~\text{odd},
\end{cases}
\qquad &&
r(x_k)=
\begin{cases}
B(x_k), & k~\text{even},
\\
A(x_k), & k~\text{odd},
\end{cases} &
\end{alignat*}
and $-m(x_k)=u(x_k)+v(x_k)+t(x_k)+r(x_k)$.
A similar relation can be found for any value of $\alpha$.
Moreover, it is possible to obtain another 5-term dif\/ference equation by considering the
alternative grid
\begin{gather*}
\widetilde{x}_{k}=(-1)^{k}(h-k/2-1/4)-1/4,
\qquad
k\in\mathbb{Z},
\end{gather*}
and proceeding along the same lines.

\section{The CBI algebra}\label{section5}

The Bannai--Ito polynomials have as an underlying algebraic structure the so-called BI
algebra~\cite{Vinet-2012}, which corresponds to a~$q\rightarrow-1$ limit of the Askey--Wilson
($AW(3)$) algebra~\cite{Zhedanov-1991}.
The algebra $AW(3)$ and the related concept of Leonard
pairs~\cite{Terwilliger-2001,Terwilliger-2008,Vidunas-2008}, describe polynomials which possess the
Leonard duality.
In this section, we obtain the algebraic structure that encodes the properties of the CBI
polynomials.

We begin by a~formal def\/inition of the CBI algebra.
\begin{definition}
The complementary Bannai--Ito (CBI) algebra is generated by the elements $\kappa_1$, $\kappa_2$,
$\kappa_3$ and the involution $r$ satisfying the relations
\begin{gather}
[\kappa_1,r]=0,
\qquad
\{\kappa_2,r\}=2\delta_3,
\qquad
\{\kappa_3,r\}=0,
\qquad
[\kappa_1,\kappa_2]=\kappa_3,
\nonumber
\\
[\kappa_1,\kappa_3]=
\frac{1}{2}\{\kappa_1,\kappa_2\}-\delta_2 \kappa_3r-\delta_3 \kappa_1r+\delta_1 \kappa_2-\delta_1\delta_3 r,
\qquad
r^2=\mathbb{I},
\nonumber
\\
[\kappa_3,\kappa_2]=
\frac{1}{2}\kappa_2^2+\delta_2 \kappa_2^2r+2\delta_3 \kappa_1r+2\delta_3 \kappa_3r+\kappa_1+\delta_4 r+\delta_5,
\label{CBI-Algebra}
\end{gather}
where $[x,y]=xy-yx$ and $\{x,y\}=xy+yx$.
The CBI algebra~\eqref{CBI-Algebra} admits the Casimir operator
\begin{gather}
Q=\frac{1}{2}\{\kappa_2^2,\kappa_1\}-\frac{\delta_2}{2} \kappa_2^2r+\kappa_1^2-\kappa_3^2+(\delta_1-1/4) \kappa_2^2
\nonumber
\\
\phantom{Q=}
{}+(\delta_3-\delta_2) \kappa_1r+2\delta_5 \kappa_1+(\delta_1\delta_3-\delta_2\delta_5) r,
\label{CBI-Casimir}
\end{gather}
which commutes with all the generators.
\end{definition}

We def\/ine the operators
\begin{gather}
\label{realization}
K_1=\mathcal{D}_{\alpha},
\qquad
K_2=x,
\end{gather}
where $K_2$ is the operator multiplication by $x$ and where $\mathcal{D}_{\alpha}$ is as given
by~\eqref{Full--OP}.
We introduce the involution~\cite{Genest-2012-1}
\begin{gather}
\label{realization-2}
P=R+\frac{\rho_2}{x}(\mathbb{I}-R).
\end{gather}
It is easily seen that $P^2=\mathbb{I}$.
Finally, we def\/ine a~fourth operator $K_3$ as follows:
\begin{gather}
\label{realization-3}
K_3=A(x)T^{+}-B(x)T^{-}+[\alpha(x-\rho_2)-2x C(x)]R-(1+2x)D(x)T^{+}R.
\end{gather}
A direct computation shows that the operators $K_1$, $K_2$ and $K_3$, together with the involution
$P$, realize the CBI algebra~\eqref{CBI-Algebra} under the identif\/ications
\begin{gather*}
K_1=\kappa_1,
\qquad
K_2=\kappa_2,
\qquad
K_3=\kappa_3,
\qquad
P=r.
\end{gather*}
The structure constants take the form
\begin{gather}
\delta_{1}=\alpha(g-\alpha+1),
\qquad
\delta_2=g-2\alpha+3/2,
\qquad
\delta_3=\rho_2,
\qquad
\delta_5=\alpha(\rho_2-1/2)+\omega/8,\nonumber
\\
\delta_4=\alpha\big(2\rho_2^2-\rho_2+1/2\big)+\rho_2 \omega/4+(8\rho_1r_1r_2+4r_1r_2-2\rho_1+2r_1+2r_2-3)/8.
\label{struct}
\end{gather}
It is worth pointing out that even though the BI and CBI polynomials can be obtained from one
another by a~Christof\/fel (resp.\ Geronimus) transformation and that they can both be obtained from
the Askey--Wilson polynomials by very similar $q\rightarrow-1$ limits, their underlying algebraic
structure are very dissimilar~\cite{Vinet-2012}.
In the realization~\eqref{realization}--\eqref{realization-3} the Casimir
operator~\eqref{CBI-Casimir} acts a~multiple of the identity{\samepage
\begin{gather*}
Q f(x)=q f(x),
\end{gather*}
where $q$ is a~complicated function of the f\/ive parameters $\rho_1$, $\rho_2$, $r_1$, $r_2$ and $\alpha$.}

The realization~\eqref{realization}--\eqref{realization-3} can be used to
obtain irreducible representations of the algebra~\eqref{CBI-Algebra} in two ``dual'' bases.
In the f\/irst basis $\{v_n,\,n\in\mathbb{N}\}$, the operator $\kappa_1$ is diagonal:
\begin{gather*}
\kappa_1v_n=\Lambda_{n}^{(\alpha)}v_n,
\end{gather*}
where $\Lambda_{n}^{(\alpha)}$ is given by~\eqref{Eigenvalues}.
Since $\kappa_1$ and $r$ commute, the operator $r$ can also be taken diagonal in this
representation.
Since $r^2=\mathbb{I}$, one f\/inds
\begin{gather*}
rv_n=\epsilon(-1)^{n}v_{n},
\end{gather*}
where $\epsilon=\pm1$ is a~representation parameter.
Given the fact that the representation parame\-ter~$\epsilon$ is only a~global multiplication factor
of~$r$, one can choose $\epsilon=1$ without loss of generality.
Because~$r$ is diagonal in the basis~$v_n$, the matrix elements of~$\kappa_2$ in the basis~$v_n$
can be calculated in a~way similar to the one employed to obtain the representations of the Hahn
algebra~\cite{Zhedanov-1992}, with additional parity requirements.
It is straightforward to show that in the basis~$v_n$, upon choosing the initial condition $a_0=0$,
the operator~$\kappa_2$ is tridiagonal with the action
\begin{gather*}
\kappa_2v_n=a_{n+1}v_{n+1}+b_nv_n+a_nv_{n-1},
\end{gather*}
where we have
\begin{gather}
\label{Matrix-Elem}
a_{n}=\sqrt{\tau_n},
\qquad
b_n=(-1)^{n}\rho_2,
\end{gather}
with $\tau_n$ given as in~\eqref{Recurrence--CBI}.
We thus have the following result.
\begin{proposition}
Let $V$ be the infinite-dimensional $\mathbb{C}$-vector space spanned by the basis vectors
$\{v_{n}\rvert n\in\mathbb{N}\}$ endowed with the actions
\begin{gather*}
\kappa_1v_{n}=\Lambda_{n}^{(\alpha)}v_{n},
\qquad
rv_{n}=(-1)^{n}v_{n},
\\
\kappa_2v_{n}=\sqrt{\tau_{n+1}} v_{n+1}+(-1)^{n}\rho_2 v_{n}+\sqrt{\tau_{n}} v_{n-1},
\\
\kappa_3v_{n}=
\big(\Lambda_{n+1}^{(\alpha)}
-\Lambda_{n}^{(\alpha)}\big)\sqrt{\tau_{n+1}}v_{n+1}
-\big(\Lambda_{n}^{(\alpha)}-\Lambda_{n-1}^{(\alpha)}\big)\sqrt{\tau_{n}}v_{n-1},
\end{gather*}
where $\Lambda_{n}^{(\alpha)}$ and $\tau_{n}$ are given by~\eqref{Eigenvalues}
and \eqref{Recurrence--CBI}, respectively.
Then $V$ is a~module for the CBI algebra~\eqref{CBI-Algebra} with structure constants taking the
values~\eqref{struct}.
The module is irreducible if none of the truncation conditions~\eqref{CBI-Trunca-Even}
and \eqref{CBI-Trunca-Odd} are satisfied.
\end{proposition}
\begin{proof}
The above considerations show that $V$ is indeed a~CBI-module.
The irreducibility stems from the fact that if the none of the truncation
conditions~\eqref{CBI-Trunca-Even} and \eqref{CBI-Trunca-Odd} are satisf\/ied, then~$\tau_{n}$ is
never zero.
\end{proof}

\begin{Corollary}
If one of the truncation conditions~\eqref{CBI-Trunca-Even} or~\eqref{CBI-Trunca-Odd} is
satisfied, then $V$ is no longer irreducible.
One can restrict to the subspace spanned by the basis vectors $\{v_{n}\,\rvert\, n=0,\ldots,N\}$
and obtain a~$(N+1)$-dimensional irreducible CBI-module.
\end{Corollary}

Thus the CBI algebra admits inf\/inite-dimensional representations where $\kappa_1$, $r$ are
diagonal and $\kappa_2$ is tridiagonal with matrix elements $\eqref{Matrix-Elem}$.
Is is readily checked that
\begin{gather*}
PI_n(x)=(-1)^nI_n(x)
\end{gather*}
and hence it is clear that the basis vectors $v_n$ correspond to the CBI polynomials themselves
\begin{gather*}
v_n=I_n(x).
\end{gather*}

Alternatively, we can consider the ``dual" basis $\{\psi_k,\,k\in\mathbb{Z}\}$, in which the
operator $\kappa_2$ is diagonal
\begin{gather*}
\kappa_2\psi_{k}=\vartheta_k\psi_k,
\end{gather*}
with the Bannai--Ito spectrum
\begin{gather}
\label{def-1}
\vartheta_k=(-1)^{k}(k/2+t+1/4)-1/4,
\end{gather}
where $t$ an arbitrary real constant.
In this basis, the involution $r$ cannot be diagonal.
Let $A_{\ell k}$ be the matrix elements of $r$ in the basis $\psi_k$.
We have
\begin{gather*}
r\psi_k=\sum_{\ell}A_{\ell k}\psi_{\ell}.
\end{gather*}
Written in the basis $\psi_k$, the anticommutation relation $\{\kappa_2,r\}=2\rho_2$ has the simple
form
\begin{gather}
\label{Dompe}
\sum_{\ell}A_{\ell,k}\{\vartheta_{\ell}+\vartheta_{k}\}\psi_{\ell}=2\rho_2\psi_{k}.
\end{gather}
For $\ell=k$, this yields
\begin{gather*}
A_{2k,2k}=\frac{\rho_2}{k+t},
\qquad
A_{2k+1,2k+1}=-\frac{\rho_2}{k+t+1}.
\end{gather*}
When $\ell\neq k$, the equation~\eqref{Dompe} reduces to
\begin{gather*}
A_{k,\ell}\{\vartheta_{k}+\vartheta_{\ell}\}=0.
\end{gather*}
From the def\/inition~\eqref{def-1} of the eigenvalues $\vartheta_{k}$, one notes that
\begin{gather}
\label{Block}
\vartheta_{2k+1}+\vartheta_{2k+2}=0.
\end{gather}
It follows from~\eqref{Block} that in the basis $\psi_k$, the operator $r$ is block diagonal with
all blocks $2\times2$.
Upon demanding that the other commutation relations of~\eqref{CBI-Algebra} be satisf\/ied, it can
be shown~\cite{Genest-2012-1} that in this basis, the operator $\kappa_1$ becomes 5-diagonal.
This result is expected since the CBI polynomials obey a~5-term dif\/ference equation of the
form~\eqref{CBI-Difference} on the Bannai--Ito grid.

We have obtained that the CBI polynomials are eigenfunctions of a~one-parameter family of operators
of the form~\eqref{Full--OP} and that two operators $\mathcal{D}_{\alpha}$, $\mathcal{D}_{\beta}$
of this family are related by the ``hidden'' symmetry operator of the CBI polynomials given
by~\eqref{Hidden-Eigen}.
In the CBI algebra, the transformation $\mathcal{D}_{\alpha}\rightarrow\mathcal{D}_{\alpha+\beta}$
is equivalent to def\/ining
\begin{gather}
\label{Transformation}
\widetilde{K}_1=K_1+\frac{\beta}{2}(\mathbb{I}-P),
\end{gather}
while leaving $K_2$ and $P$ unchanged.
The operator $K_3$ is transformed to
\begin{gather*}
\widetilde{K}_3=K_3-\beta PK_2+\beta\delta_3.
\end{gather*}
Upon using $\widetilde{K}_2=K_2$, one f\/inds that the algebra becomes
\begin{gather*}
[\widetilde{K}_1,P]=0,
\qquad
\{\widetilde{K}_2,P\}=2\widetilde{\delta}_3,
\qquad
\{\widetilde{K}_3,P\}=0,
\qquad
[\widetilde{K}_1,\widetilde{K}_2]=\widetilde{K}_3,
\\
[\widetilde{K}_1,\widetilde{K}_3]=
\frac{1}{2}\{\widetilde{K}_1,\widetilde{K}_2\}
-\widetilde{\delta}_2\widetilde{K}_3P
-\widetilde{\delta}_3\widetilde{K}_1P
+\widetilde{\delta}_1\widetilde{K}_2
-\widetilde{\delta}_1\widetilde{\delta}_3P,
\\
[\widetilde{K}_3,\widetilde{K}_2]=
\frac{1}{2}\widetilde{K}_2^2
+\widetilde{\delta}_2\widetilde{K}_2^2P
+2\widetilde{\delta}_3\widetilde{K}_1P
+2\widetilde{\delta}_3\widetilde{K}_3P
+\widetilde{K}_1+\widetilde{\delta}_4P+\widetilde{\delta}_5,
\end{gather*}
with the structures constants
\begin{gather*}
 \widetilde{\delta}_1=\delta_1+\beta(\delta_2-1/2),
\qquad
\widetilde{\delta}_{2}=\delta_2-2\beta,
\qquad
\widetilde{\delta_3}=\delta_3,
\\
 \widetilde{\delta}_4=\delta_4+\beta(2\delta_3^2-\delta_3+1/2),
\qquad
\widetilde{\delta}_5=\delta_5+\beta(\delta_3-1/2).
\end{gather*}
It is thus seen that the transformation~\eqref{Transformation} leaves the general form of the CBI
algebra~\eqref{CBI-Algebra} unaf\/fected and corresponds only to a~change in the structure
parameters.

\section{Three OPs families related to the CBI polynomials}\label{section6}

In this section, we exhibit the relationship between the complementary Bannai--Ito polynomials
and three other families of orthogonal polynomials: the recently discovered dual $-1$
Hahn~\cite{VZhedanov-2011,Z-2011} and para-Krawtchouk polynomials~\cite{VZ-2012} and the classical
symmetric Hahn polynomials.

\subsection[Dual $-1$ Hahn polynomials]{Dual $\boldsymbol{-1}$ Hahn polynomials}\label{section6.1}

The dual $-1$ Hahn polynomials have been introduced in~\cite{VZhedanov-2011} as $q=-1$ limits of
the dual $q$-Hahn polynomials.
They have appeared in the context of perfect state transfer in spin chains~\cite{Z-2011} and also
as the Clebsch--Gordan coef\/f\/icients of the $sl_{-1}(2)$ algebra
in~\cite{Genest-2012-2,Vinet-2011}.
Moreover, the $-1$ Hahn polynomials have occurred, in their symmetric form, as wavefunctions for
f\/inite parabosonic oscillator models~\cite{VDJ-2011,VDJ-2011-2}.
These polynomials\footnote{To recover the formulas found in~\cite{VZhedanov-2011},
a~re-parametrization is necessary.}, denoted $Q_n(x)$, can be obtained from the CBI polynomials
through the limit $\rho_1\rightarrow\infty$.

Taking the limit $\rho_1\rightarrow\infty$ in the~\eqref{Recurrence--CBI}, one obtains the
recurrence relation of the monic dual $-1$ Hahn polynomials
\begin{gather*}
Q_{n+1}(x)+(-1)^{n}\rho_2 Q_{n}(x)+\sigma_{n}Q_{n-1}(x)=x Q_{n}(x),
\end{gather*}
where $r_n$ has the expression
\begin{gather*}
\sigma_{2n}=-n(n-r_1-r_2),
\qquad
\sigma_{2n+1}=-(n+\rho_2-r_1+1/2)(n+\rho_2-r_2+1/2).
\end{gather*}
The polynomials $Q_n(x)$ have the hypergeometric representation
\begin{gather*}
Q_{2n}(x)=\xi_{2n}\; \pFq{3}{2}{-n,\rho_2+x,\rho_2-x}{\rho_2-r_1+1/2,\rho_2-r_2+1/2}{1},
\\
Q_{2n+1}(x)=
\xi_{2n+1} (x-\rho_2)\;\pFq{3}{2}{-n,\rho_2+x+1,\rho_2-x+1}{\rho_2-r_1+3/2,\rho_2-r_2+3/2}{1},
\end{gather*}
with normalization coef\/f\/icients
\begin{gather*}
\xi_{2n}=(\rho_2-r_1+1/2)_{n}(\rho_2-r_2+1/2)_n,
\qquad
\xi_{2n+1}=(\rho_2-r_1+3/2)_{n}(\rho_2-r_2+3/2)_n.
\end{gather*}
These formulas are obtained from~\eqref{CBI-Hypergeo} in the same limit.
Dividing~\eqref{Full--OP} by $\rho_1$ and taking the limit $\rho_1\rightarrow\infty$, one f\/inds
that the polynomials $Q_n(x)$ satisfy the eigenvalue equation
\begin{gather*}
\mathcal{E}^{(\alpha)}Q_{n}(x)=\nu_{n}^{(\alpha)}Q_n(x),
\end{gather*}
with eigenvalues
\begin{gather*}
\nu_{2n}^{(\alpha)}=n,
\qquad
\nu_{2n+1}^{(\alpha)}=n+\alpha.
\end{gather*}
The operator $\mathcal{E}^{(\alpha)}$ is found to be
\begin{gather*}
\mathcal{E}^{(\alpha)}=\mathcal{E}^{(0)}+\alpha\frac{(x-\rho_2)}{2x}(\mathbb{I}-R),
\end{gather*}
where
\begin{gather*}
\mathcal{E}^{(0)}=I(x) T^{+}+J(x) T^{-}+K(x) R+L(x) T^{+}R- (I(x)+J(x)+K(x)+L(x))\mathbb{I}.
\end{gather*}
The coef\/f\/icients are given by
\begin{gather*}
I(x)=\frac{(x+\rho_2+1)(2x-2r_1+1)(2x-2r_2+1)}{8(x+1)(2x+1)},
\\
J(x)=\frac{(\rho_2-x)(2x+2r_1-1)(2x+2r_2-1)}{8x(2x-1)},
\\
K(x)=\frac{(x-\rho_2)(4x^2+4r_1r_2-1)}{4x(4x^2-1)},
\\
L(x)=\frac{\rho_2(2x-2r_1+1)(2x-2r_2+1)}{8x(x+1)(2x+1)}.
\end{gather*}
Lastly, it is seen that in the limit $\rho_1\rightarrow\infty$, the CBI algebra becomes
\begin{gather*}
[\kappa_1,r]=0,
\qquad
\{\kappa_2,r\}=2\gamma_3,
\qquad
\{\kappa_3,r\}=0,
\qquad
[\kappa_1,\kappa_2]=\kappa_3,
\\
[\kappa_1,\kappa_3]=\gamma_1 \kappa_2-\gamma_1\gamma_3 r-\gamma_2 \kappa_3r,
\qquad
[\kappa_3,\kappa_2]=\gamma_2 \kappa_2^2 r+2\gamma_3 \kappa_1r+2\gamma_3 \kappa_3r+\kappa_1+\gamma_4 r+\gamma_5,
\end{gather*}
where we have identif\/ied $\kappa_1=\mathcal{E}^{(\alpha)}$, $\kappa_2=x$ and $P=r$ with $P$ given
by~\eqref{realization-2}.
The structure parameters have the expression
\begin{gather*}
\gamma_1=\alpha(1-\alpha),
\qquad
\gamma_2=1-2\alpha,
\qquad
\gamma_3=\rho_2,
\\
\gamma_4=\alpha\big(2\rho_2^2-\rho_2+1/2\big)+\rho_2(1-r_2-r_1)+r_1r_2-1/4,
\\
\gamma_5=(2\alpha\rho_2-\alpha-r_1-r_2+1)/2.
\end{gather*}
Other properties of the polynomials $Q_{n}(x)$ can be obtained directly using the limiting
procedure.

\subsection{The symmetric Hahn polynomials}\label{section6.2}
It is possible to relate the CBI polynomials to the symmetric Hahn
polynomials~\cite{Karlin-1961,Koekoek-2010} through a~direct identif\/ication of the CBI parameters.
This identif\/ication can be performed in three dif\/ferent ways by examining the cases for which
the def\/ining operator $\mathcal{D}_{\alpha}$~\eqref{Full--OP} of the CBI polynomials reduces to
a~classical three-term dif\/ference operator involving only the discrete shifts~$T^{+}$,~$T^{-}$
and the identity operator $\mathbb{I}$.

We consider the operator $\mathcal{D}_{\alpha}$ in~\eqref{Full--OP} with the following parameter
identif\/ication:
\begin{gather}
\label{Pp1}
\rho_1=-\frac{1}{2},
\qquad
\rho_2=0,
\qquad
\alpha=\frac{1}{2}(1-r_1-r_2).
\end{gather}
With these values of the parameters, the eigenvalue equation~\eqref{Eigen-Equation} reduces to
\begin{gather}
\label{Reduction}
B(x)I_n(x+1)-(B(x)+D(x))I_n(x)+D(x)I_{n}(x-1)=\lambda_n I_{n}(x),
\end{gather}
with coef\/f\/icients
\begin{gather*}
B(x)=(x-r_1+1/2)(x-r_2+1/2),
\qquad
D(x)=(x+r_1-1/2)(x+r_2-1/2),
\end{gather*}
and eigenvalues
\begin{gather*}
\lambda_n=n(n-2r_1-2r_2+1).
\end{gather*}
We consider the parametrization
\begin{gather}
\label{Para}
r_1=\frac{N+1}{2},
\qquad
\alpha^*=\beta^{*}=-r_1-r_2,
\end{gather}
where $N$ is an even or odd integer.
It is seen from~\eqref{CBI-Trunca-Even} and \eqref{CBI-Trunca-Odd} that~\eqref{Para} is an
admissible truncation condition for both parities of $N$.
Upon introducing the variable $\widetilde{x}=x+r_1-1/2$, the coef\/f\/icients of the eigenvalue
equation~\eqref{Reduction} become
\begin{gather*}
B(x)=(\widetilde{x}-N)\big(\widetilde{x}+\alpha^*+1\big),
\qquad
D(x)=\widetilde{x}\big(\widetilde{x}-\beta^*-N-1\big),
\end{gather*}
and the eigenvalues have the expression
\begin{gather*}
\lambda_n=n\big(n+\alpha^*+\beta^*+1\big).
\end{gather*}
This corresponds to the dif\/ference equation of the Hahn polynomials~\cite{Koekoek-2010}.
With the parametrization~\eqref{Pp1} and \eqref{Para}, the recurrence relation~\eqref{Recurrence-1}
of the CBI polynomials becomes
\begin{gather*}
I_{n+1}(x)+\omega_n I_{n-1}(x)=xI_{n}(x),
\end{gather*}
where
\begin{gather*}
\omega_n=
\frac{n(N-n+1)(n+\alpha^*+\beta^*)(n+\alpha^{*}+\beta^{*}+N+1)}
{4(2n+\alpha^{*}+\beta^{*}-1)(2n+\alpha^{*}+\beta^{*}+1)},
\end{gather*}
which is indeed the recurrence relation of symmetric Hahn polynomials.
A simple calculation shows that upon taking the parametrization~\eqref{Para} in structure
parameters~\eqref{struct}, one has
\begin{gather*}
\delta_2=\delta_3=\delta_4=0,
\end{gather*}
and hence the algebra reduces to
\begin{gather*}
[K_1,P]=0,
\qquad
\{K_2,P\}=0,
\qquad
\{K_3,P\}=0,
\qquad
[K_1,K_2]=K_3,
\\
[K_1,K_3]=\frac{1}{2}\{K_1,K_2\}+\delta_1K_2,
\qquad
[K_3,K_2]=\frac{1}{2}K_2^2+K_1+\delta_5,
\end{gather*}
with the remaining structure parameters
\begin{gather*}
\delta_1=\frac{1}{4}(r_1+r_2),
\qquad
\delta_5=\frac{1}{4}(r_1-1/2)(r_2-1/2).
\end{gather*}
Thus we recover the Hahn algebra since the involution $P$ no longer plays a~determining role.
For reference, we record the two following alternate choices of CBI parameters which also lead to
symmetric Hahn polynomials
\begin{gather*}
\rho_2=r_1=0,
\qquad
\alpha=\frac{1}{4}(2\rho_1-2r_2+3),
\qquad
\text{or}
\qquad
\rho_2=r_2=0,
\qquad
\alpha=\frac{1}{4}(2\rho_1-2r_1+3).
\end{gather*}

\subsection{Para-Krawtchouk polynomials}\label{section6.3}
The para-Krawtchouk polynomials have been found~\cite{VZ-2012} in the design of spin chains
ef\/fecting perfect quantum state transfer.
These polynomials are directly connected to the complementary Bannai--Ito polynomials through the
identif\/ication
\begin{gather*}
\rho_1=\frac{\gamma-N-3}{4},
\qquad
\rho_2=0,
\qquad
r_1=\frac{N+1+\gamma}{4},
\qquad
r_2=0,
\end{gather*}
when $N$ is a~positive odd integer.
(When $N$ is an even integer, the para-Krawtchouk are directly related to the Bannai--Ito
polynomials.) In~\cite{VZ-2012}, the eigenvalue equation for the para-Krawtchouk polynomials was
found; this operator corresponds to the operator~\eqref{Full--OP} with a~specif\/ic value of the
free parameter
\begin{gather*}
\alpha=\frac{1-N}{4}.
\end{gather*}
It is interesting to note that in the case $\rho_2=0$, the CBI polynomials and their descendants
become symmetric and thus the ``hidden'' eigenvalue equation~\eqref{Hidden-Eigen} appears trivial.
Notwithstan\-ding this, the corresponding symmetric polynomials are still eigenfunctions of a~Dunkl
operator with a~free parameter in the odd sector of the spectrum.

\section{Conclusion}\label{section7}
We have presented a~systematic study of the complementary Bannai--Ito polynomials.
We showed that these OPs are eigenfunctions of a~one-parameter family of second order Dunkl shift
operators and that in consequence they satisfy a~one-parameter family of f\/ive-term dif\/ference
equations on grids of the Bannai--Ito type.
This result makes explicit the bispectrality of the CBI polynomials and places this OPs family
outside the scope of the Leonard duality.
Moreover, we have obtained the algebraic structure associated to the CBI polynomials which we named
the complementary Bannai--Ito algebra.
It was observed that this quadratic algebra is a~deformation of the Askey--Wilson algebra with an
involution.
Lastly, we identif\/ied how the CBI polynomials are related to three other families of OPs.

The investigation of the continuum limit of the BI polynomials has led to connections with other
families of $-1$ orthogonal polynomials which satisfy f\/irst order dif\/ferential/Dunkl equations.
It is hence of interest to examine the continuum limit of the CBI polynomials; this question will
be treated in a~future publication.

\subsection*{Acknowledgements}
V.X.G.~holds a~scholarship from \emph{Fonds de recherche qu\'eb\'ecois~-- nature et technologies}
(FRQNT).
The research of L.V.\ is supported in part by the \emph{Natural Science and Engineering Council of
Canada} (NSERC).
A.Z.~would like to thank the \emph{Centre de Recherches Math\'ematiques} (CRM) for its hospitality.

\pdfbookmark[1]{References}{ref}
\LastPageEnding

\end{document}